\documentclass{amsart}

\usepackage{mathrsfs}
\usepackage{amssymb}
\usepackage[usenames,dvipsnames]{xcolor}
\usepackage{mathtools}
\usepackage{caption}
\usepackage{enumerate}
\usepackage{natbib}
\usepackage{hyperref}
\usepackage[nameinlink, capitalize, noabbrev]{cleveref}
\usepackage{stmaryrd}
\usepackage{tikz}
\usetikzlibrary{intersections}
\usepackage{pgfplots}
\pgfplotsset{width=10cm,compat=1.9}
\usepgfplotslibrary{external}

\hypersetup{ colorlinks=true,linkcolor=teal,    citecolor=magenta, }

\usepackage{comment}

\def\set#1{{\def\st{\;:\;}\left\{#1\right\}}}
\def\abs#1{\left\vert{#1}\right\vert}
\def \<#1>{{\left\langle{#1}\right\rangle}}

\def\PP{\mathbb P}
\def\OK{\mathcal{O}_K}
\def\FF{\mathbb F}
\def\CC{\mathbb C}
\def\NN{\mathbb N}
\def\RR{\mathbb R}
\def\Q-{\overline{\mathbb Q}}

\DeclareMathOperator{\Aut}{Aut}
\DeclareMathOperator{\FPP}{FPP}
\DeclareMathOperator{\Per}{Per}

\DeclareMathOperator{\IMG}{IMG}

\DeclareMathOperator{\St}{St}
\DeclareMathOperator{\Sym}{Sym}

\DeclareMathOperator{\st}{st}

\newtheorem{Theorem}{Theorem}

\newtheorem{proposition}{Proposition}[section]
\newtheorem{lemma}[proposition]{Lemma}

\newtheorem{theorem}[proposition]{Theorem}

\newtheorem{Question}[proposition]{Question}
\theoremstyle{definition}

\title{On the fixed-point proportion of self-similar groups}
\author{Jorge Fariña-Asategui and Santiago Radi}
\address{Jorge Fariña-Asategui: Centre for Mathematical Sciences, Lund University, 223 62 Lund, Sweden -- Department of Mathematics, University of the Basque Country UPV/EHU, 48080 Bilbao, Spain}
\email{jorge.farina\_asategui@math.lu.se}
\address{Santiago Radi: Department of Mathematics, Texas A\&M University, 77843 College Station, U.S.A.
}
\email{santiradi@tamu.edu}
\thanks{The first author is supported by the Spanish Government, grant PID2020-117281GB-I00, partly with FEDER funds. The author also acknowledges support from the Walter Gyllenberg Foundation from the Royal Physiographic Society of Lund. The second author is supported by Grigorchuk's Simons Foundation Grant MP-TSM-00002045 and the department of Mathematics of Texas A\&M University.}

\begin{document}

\begin{abstract}
We prove that super strongly fractal groups acting on regular rooted trees have null fixed-point proportion. In particular, we show that the fixed-point proportion of an infinite family of iterated monodromy groups of exceptional complex polynomials have the same property. The proof uses the approach of Rafe Jones in \cite{jones_arborealsurvey} based on martingales and a recent result of the first author on the dynamics of self-similar groups \cite{JorgeCyclicity}.
\end{abstract}

\maketitle

\section{Introduction}
\label{section: introduction}

Given an integer $a_0\in \mathbb{Z}$ and a polynomial $f\in \mathbb{Z}[x]$, one may define a sequence $\{f^n(a_0)\}_{n\ge 0}$ by considering iterates of the polynomial $f$ on $a_0$. Many interesting families of numbers arise in this manner, such as Fermat numbers and Mersenne numbers \cite{jones_arborealsurvey}. One of the main longstanding open problems in arithmetic dynamics concerns which of the terms of such a sequence $\{f^n(a_0)\}_{n\ge 0}$ are prime. A related easier problem is understanding the density of prime divisors of the terms in the sequence  $\{f^n(a_0)\}_{n\ge 0}$. More precisely, if $P_f(a_0)$ denotes the set of primes dividing some non-zero term of the sequence $\{f^n(a_0)\}_{n\ge 0}$, one is interested in computing the natural density:
$$\mathcal{N}(P_f(a_0)):=\lim_{k\to \infty}\frac{\# \{p\in P_f(a_0): p\le k\}}{\#\{p\text{ prime}: p\le k\}}.$$

The pioneering work of Odoni in \cite{Odoni_1985}, settled an approach to give an upper bound for the quantity $\mathcal{N}(P_f(a_0))$. This approach is to consider the Galois group $G_\infty(f)$ associated to the preimages of $0$ by all the iterates of the polynomial $f$ and then take the natural action of $G_\infty(f)$ on the tree of preimages of $0$ via the iterates of~$f$. Using the \v{C}ebotarev Density Theorem, one can give an upper bound of $\mathcal{N}(P_f(a_0))$ via the so-called \textit{fixed-point proportion} of $G_\infty(f)$, i.e. the proportion of elements in $G_\infty(f)$ fixing at least one infinite path in this tree of preimages. We use the notation $\FPP(G)$ to denote the fixed-point proportion of a group $G$.

Even though the original work of Odoni in \cite{Odoni_1985} concerned the specific function $f(x) = x^2-x+1$ and the point $a_0=2$, his work was later generalized by Jones to any rational function and any rational point. Furthermore, Jones and Manes extended this approach to any number field $K$ and any rational function $f\in K(x)$: 
\begin{theorem}[{see {\cite[Theorem 6.1]{JonesManes2012}}}]
\label{theorem: Jones sequences}
Let $K$ be any number field, $M_K^0$ the set of prime ideals in $\OK$, $v_\mathfrak{p}$ the valuation in $\OK$ for $\mathfrak{p}\in M_K^0$, $f$ any rational function in $K(x)$ of degree at least $2$ and $a_0 \in K$. Denote by $f^n$ the $n$-th iteration of $f$ and define the set 
$$P_f(a_0) := \set{\mathfrak{p} \in M_K^0: v_\mathfrak{p}(f^n(a_0)) > 0 \text{ for some $n$ such that $f^n(a_0) \neq 0$}};$$
namely, the set of prime ideals that divide at least one non-zero term in the orbit of $a_0$ by $f$. If $0$ does not appear infinitely many times in the orbit of $\infty$, then $$\mathcal{D}(P_f(a_0)) \leq \FPP(G_\infty(K,f,0)),$$ where $\mathcal{D}$ is the Dirichlet density in $K$ and $G_\infty(K,f,0)$ is a Galois group associated to $f$ and $K$.
\end{theorem}

It has been shown that the fixed point proportion of groups acting on regular rooted trees has further applications to arithmetic dynamics other than \cref{theorem: Jones sequences}. For instance, there is a link between the number of periodic points of a rational function $f$ over finite fields and the fixed-point proportion of the iterated monodromy group of $f$:

\begin{theorem}[{see {\cite[Theorem 3.11]{BridyJones2022}, \cite[Proposition 5.3]{juul2014wreath} and \cite[Proposition~6.4.2]{Self_similar_groups}}}]
\label{theorem: jones periodic points}
Let $K$ be a number field and $f$ a rational function in $K(x)$ of degree at least $2$. For any prime ideal $\mathfrak{p}$ in $\OK$ denote $f_\mathfrak{p}$ the reduction of $f$ to $\FF_\mathfrak{p} := \OK/\mathfrak{p}$. Denote $\Per(f_\mathfrak{p},\FF_\mathfrak{p})$ the set of periodic points of $f$ in the finite field~$\FF_\mathfrak{p}$. If $f$ is post-critically finite, then 
\begin{equation*}
\liminf_{N(\mathfrak{p}) \rightarrow \infty} \frac{\# \Per(f,\FF_\mathfrak{p})}{N(\mathfrak{p})+1} \leq \FPP(\IMG(f))
\end{equation*}
where $\IMG(f)$ is the iterated monodromy group of $f$.
\end{theorem}

Further applications can be found in \cite[Corollary 3.5]{BridyJones2022} and \cite[Theorem 1.4]{Jones2007}. 

The above results motivate the study of the fixed-point proportion of groups acting on regular rooted trees. In view of the upper bound in \cref{theorem: Jones sequences} and \cref{theorem: jones periodic points}, it is of special interest when the fixed-point proportion is zero:

\begin{Question}
\label{question_FPPgroups}
Let $T$ be a regular rooted tree and $G$ a subgroup of $\Aut(T)$:
\begin{enumerate}[\normalfont(i)]
\item What can be said about the value of $\FPP(G)$?
\item  What assumptions must the group $G$ satisfy so that $\FPP(G) = 0$?
\end{enumerate}
\end{Question}

Regarding \cref{question_FPPgroups}\textcolor{teal}{(ii)} it is known that iterated wreath products have null fixed-point proportion if and only if they are level-transitive; see \cite[Corollary 2.6]{AbertVirag2004} and \cite[Corollary 4]{Radi2024}. Moreover, Jones in \cite[Theorem 1.5]{jones2012fixedpointfree} provided sufficient conditions for contracting groups to have null fixed-point proportion. However, there is a gap in the proof; see \cite{Corrigendum}. Our main result generalizes the case of iterated wreath products to a wider family of groups, known as super strongly fractal groups; see \cref{section: preliminaries} for the unexplained terms here and elsewhere in the introduction.

\begin{Theorem}
\label{theorem: FPP SSF}
Let $G \leq \Aut(T)$ be a super strongly fractal group. Then $$\FPP(G) = 0.$$
\end{Theorem}

The proof of \cref{theorem: FPP SSF} follows the same strategy based on martingales as in \cite{BridyJones2022,jones2012fixedpointfree}. However, for the key step instead of using the contracting properties of the group, we use a completely new approach based on the ergodic properties of self-similar groups developed by the first author in \cite{JorgeCyclicity}.

As it was observed in \cref{theorem: jones periodic points}, the fixed-point proportion of iterated monodromy groups is of special interest for applications to arithmetic dynamics.

The concept of dynamically exceptional was introduced by Makarov and Smirnov in \cite{MakarovSmirnov1996}. A polynomial $f$ is dynamically exceptional if there exists a non-empty set $E \subseteq \PP^1(\CC)$ such that $f ^{-1}(E) \setminus C_f = E$, where $C_f$ is the set of critical points of~$f$. It can be proved that $E$ has at most two points. When $\#E = 2$, Jones proved in \cite[Proposition 2]{jones2012fixedpointfree} that the polynomial $f$ is linearly conjugated to a Chebyshev polynomial and therefore, that the fixed-point proportion of $\IMG(f)$ is $1/2$ when $\deg(f)$ is odd and $1/4$ when $\deg(f)$ is even.

If $\deg(f)=2$, then either $f$ is not dynamically exceptional or $\#E=2$. Therefore, the dynamically exceptional case with $\deg(f)=2$ is completely settled by the work of Jones; see \cite[Corollary~1.3]{jones2012fixedpointfree}. However, for a dynamically exceptional polynomial such that $\deg(f)\ge 3$ and $\#E=1$, nothing is known. It was conjectured by the authors in \cite[page~442]{BridyJones2022} that the iterated monodromy groups of dynamically exceptional polynomials non-conjugated to Chebyshev polynomials have null fixed-point proportion. However, their methods do not suffice to compute the fixed-point proportion of any of these groups. As an application of \cref{theorem: FPP SSF}, we prove the following:

\begin{Theorem}
\label{theorem: dynamically exceptional FPP=0}
For each $d \geq 3$, there exists a dynamically exceptional polynomial $f$ with $\# E = 1$ and degree $d$ such that $$\FPP(\IMG(f)) = 0.$$
\end{Theorem}

\cref{theorem: dynamically exceptional FPP=0} provides infinitely many examples supporting this conjecture in \cite{BridyJones2022}. The proof of \cref{theorem: dynamically exceptional FPP=0} is constructive. Explicit examples of such polynomials are given at the end of \cref{sec_FPP_dyn_exc_poly}.

\subsection*{Organization}

In \cref{section: preliminaries} the necessary background is introduced. \cref{section: FPP for SSF} is devoted to prove \cref{theorem: FPP SSF}. In \cref{sec_FPP_dyn_exc_poly}, we prove \cref{theorem: dynamically exceptional FPP=0} and finally in \cref{section: Scope of theorem}, we give examples and discuss the scope of \cref{theorem: FPP SSF}.

\subsection*{Acknowledgements}

We would like to thank Volodymyr Nekrashevych for helpful discussions on iterated monodromy groups. We also like to thank Rostislav Grigorchuk and Texas A\&M for their warm hospitality while this work was being carried out. The second author would like to thank Thomas Tucker for suggesting the problem and for helpful discussions. The first author would like to thank the University of Southampton and the University of Salerno for their invitations to talk on this paper at their seminars, while the second author thanks Louisiana State University for its invitation to talk in the Southern regional number theory conference. Finally, we would like to thank the anonymous referees whose feedback contributed greatly to the final version of the paper. 

\section{Preliminaries}
\label{section: preliminaries}

\subsection{Groups acting on rooted trees}
\label{subsection: Groups acting on rooted trees}

Let $T$ be a \textit{$d$-regular rooted tree}, i.e. the infinite tree with root $\emptyset$ and whose vertices all have exactly $d$ immediate descendants. The set of vertices at a distance exactly $n \geq 1$ from the root form the \textit{$n$-th level of $T$} and will be denoted $\mathcal{L}_n$. The vertices whose distance is at most $n$ from the root form the \textit{$n$-th truncated tree}~$T^n$.

We may describe the vertices of $T$ as finite words on $X := \set{1,\dots,d}$. Using this identification, we denote by $\partial T$ the set of infinite words in $X$. The set $\partial T$ is called the \textit{boundary} of the tree $T$. The group of automorphisms of $T$, denoted $\Aut(T)$, is the group of bijective functions from $T$ to itself preserving adjacency between vertices. In particular $\Aut(T)$ acts on $\mathcal{L}_n$ for all $n \ge 1$. For any vertex $v\in T$, the subtree $T_v$ rooted at $v$ is isomorphic to $T$.

Given a vertex $v \in T$, we write $\st(v)$ for the \textit{stabilizer of the vertex $v$}, namely, the subgroup of the elements $g \in \Aut(T)$ such that $g(v) = v$. Given $n \ge 1$, we write $\St(n) := \bigcap_{v\in \mathcal{L}_n}\st(v)$, and we call it the \textit{stabilizer of level $n$}. The subgroup $\mathrm{St}(n)$ is a normal subgroup of finite index in $\mathrm{Aut}(T)$.

We can make $\Aut(T)$ a topological group by declaring $\{\St(n)\}_{n\ge 1}$ to be a base of neighborhoods of the identity. This topology is called the \textit{congruence topology}. The group $\Aut(T)$ is a profinite group with respect to the congruence topology.

Let $v \in T$, $1\le n \le \infty$ and $g \in \Aut(T)$. The unique automorphism $g|_v^n \in \Aut(T^n)$ such that for all $w \in T^n$
$$g(vw)=g(v)g|_v^n(w),$$ is called the \textit{section of $g$ at vertex $v$ of depth $n$}. For $n=1$, the permutation $g|_v^1$ is known as the \textit{label of $g$ at $v$} and for $n=\infty$ we shall drop the superscript and simply write $g|_v$. For every $n\ge 1$, we define the map $\pi_n:\Aut(T)\to \Aut(T^n)$ given by $\pi_n(g):=g|_\emptyset^n $. Sections allow us to define the following isomorphism for any natural number $n\ge 1$: 
\begin{align*}
\begin{split}
\psi_n:\Aut(T) &\to \big(\Aut(T) \times \dotsb \times \Aut(T)\big) \rtimes \Aut(T^n) \\ g &\mapsto (g|_{v_1},\dotsc,g|_{v_{d^n}}) \pi_n(g),
\end{split}
\label{equation: Aut and semidirect product}
\end{align*}
where $v_1,\dotsc,v_{d^n}$ are all the distinct vertices in $\mathcal{L}_n$ labeled from left to right.

To describe the elements of $\Aut(T)$, we will use the isomorphism $\psi_1$ and write $g=(g_1,\dotsc,g_d)\pi_1(g)$ by a slight abuse of notation.

If $v \in \mathcal{L}_n$, we define the map $\varphi_v: \Aut(T) \rightarrow \Aut(T)$ via $g \mapsto g|_v$. Note that the restriction of $\varphi_v$ to $\st(v)$ is a continuous group homomorphism. Continuity is immediate from $\varphi_v^{-1}(\St(m)) \geq \St(m+n)$ for all $m \ge 1$. 

Let us fix a subgroup $G \le \Aut(T)$. We define vertex stabilizers and level stabilizers as $\st_G(v) := \st(v) \cap G$ and $\St_G(n) := \St(n) \cap G$ for $v \in T$ and $n \ge 1$ respectively. We say that a group $G \leq \Aut(T)$ is \textit{level-transitive} if the action of $G$ on each level $\mathcal{L}_n$ is transitive for all $n \ge 1$. 

\subsection{Self-similarity and fractality}

A group $G \leq \Aut(T)$ is \textit{self-similar} if for all $v \in T$ and $g \in G$, we get $g|_v \in G$. We say that a group $G \leq \Aut(T)$ is \textit{fractal} if $G$ is self-similar, level-transitive and $\varphi_v(\st_G(v)) = G$ for all $v \in T$ and \textit{super strongly fractal} if $G$ is self-similar, level-transitive and $\varphi_v(\St_G(n)) = G$ for all $v \in \mathcal{L}_n$ and $n \ge 1$. Note that, as $\varphi_v$ is continuous for every $v\in T$, if $G$ is super strongly fractal, then its closure $\overline{G}$ in $\mathrm{Aut}(T)$ is also super strongly fractal.

Fractal properties of self-similar groups have been recently linked to the dynamics and the ergodic properties of self-similar profinite groups by the first author in \cite{JorgeCyclicity}. Furthermore, certain fractal groups have been recently used to construct scale groups; see \cite{Liftability} and the references herein.

\subsection{Fixed-point proportion}

Let $G \leq \Aut(T)$. The fixed-point proportion of $G$ is defined as 
\begin{equation*}
\FPP(G) := \lim_{n \rightarrow \infty} \frac{\# \set{g \in \pi_n(G): \text{$g$ fixes a vertex in $\mathcal{L}_n$}}}{\abs{\pi_n(G)}}.
\label{equation: FPP definition}
\end{equation*}

In \cite[Lemma 2.9]{Radi2024} it is proved that the limit in the definition of $\FPP(G)$ always exists. Moreover, note that the definition of fixed-point proportion only depends on the finite quotients $\pi_n(G)$ and thus
$$\FPP(G) = \FPP(\overline{G}),$$
as $\pi_n(G) = \pi_n(\overline{G})$ for all $n \ge 1$ by \cite[Corollary 1.1.8(c)]{Ribes2000}.

Now $\Aut(T)$ is profinite for the congruence topology, so $\overline{G}$ is also profinite. Then, we can associate to $\overline{G}$ a unique Haar probability measure $\mu$, turning $\overline{G}$ into a probability space. A well-known property of the Haar measure \cite[Lemma 16.2]{Fried1986} is that if $S$ is a measurable subset of $\overline{G}$, then 
\begin{equation*}
\mu(S) = \lim_{n \rightarrow \infty} \frac{\#\pi_n(S)}{\abs{\pi_n(\overline{G})}}.  
\label{equation: Haar measure limit}  
\end{equation*}
Then, it is immediate to see that the fixed-point proportion of a group $G$ is given by
\begin{equation*}
\FPP(\overline{G}) = \mu(\set{g \in G: g \text{ fixes a vertex in $\mathcal{L}_n$ for all $n\ge 1$}}),
\label{ec_FPP_setS}
\end{equation*}

\section{Fixed-point proportion for super strongly fractal groups}
\label{section: FPP for SSF}

In this section, we prove that super strongly fractal groups have zero fixed-point proportion. We follow the approach of Jones in \cite[page 22]{jones_arborealsurvey} based on martingales. Since a group and its closure have the same fixed-point proportion, we may prove \cref{theorem: FPP SSF} for closed subgroups of $\Aut(T)$, so that these groups are endowed with a Haar probability measure. 

\subsection{Fixed-point process and martingales}

Let $G\le \Aut(T)$ be a closed subgroup. Then for every $n\ge 1$, we define the real measurable function $X_n: G \to \NN\cup\{0\}$ given by $$X_n(g) := \#\set{v \in \mathcal{L}_n: g(v) = v}.$$
In other words, for $n \geq 1$ and $g\in G$, the number $X_n(g)$ is no more than how many vertices are fixed by $g$ at level $n$. These functions define a real stochastic process $\{X_n\}_{n \geq 1}$ called the \textit{fixed-point process of $G$}.

Among real stochastic processes in a general probability space $(X,\mu)$, we are interested in martingales. A \textit{martingale} is a real stochastic process $\set{Y_n}_{n \geq 1}$ in $(X,\mu)$ such that for all $n\ge 1$ we have:
\begin{enumerate}[\normalfont(i)]
\item $\mathbb{E}(Y_n) < \infty$ for every $n\ge1$;
\item $\mathbb{E}(Y_n \mid Y_1 = t_1,\dotsc,Y_{n-1} = t_{n-1}) = t_{n-1}$ for every $t_1,\dotsc,t_{n-1} \in \RR$ such that $\mu(Y_1 = t_1,\dotsc, Y_{n-1} = t_{n-1}) > 0$.
\end{enumerate}

The main reason for us to be interested in martingales is the following classical convergence theorem:

\begin{theorem}[{see {\cite[Theorem 7, chapter 12]{Grimmett_Random_processes}}}]
If $\set{Y_n}_{n \geq 1}$ is a non-negative martingale with $\mathbb{E}(Y_1) < \infty$, then the limit
$$\lim_{n \rightarrow + \infty} Y_n(g)$$
exists almost surely with finite value.
\label{teo_convergence_theorem_martingale}
\end{theorem}

If the fixed-point process $\set{X_n}_{n\ge 1}$ is a martingale, it satisfies the assumptions of \cref{teo_convergence_theorem_martingale} and $\{X_n\}_{n\ge 1}$ is eventually constant for almost all $g \in G$, as $X_n(g)$ is a non-negative integer for every $n\ge 1$ and $g\in G$. However, it could happen that for different non-null sets of~$G$, the constants to which their fixed-point process converge are different; see \cite{Radi2024}. In order to prove that $\{X_n\}_{n\ge 1}$ is eventually zero almost everywhere, Jones proposes a strategy in \cite{jones_arborealsurvey}. A proof of a similar strategy in a simpler case can be found in \cite{Jones2007}. For completeness, we give a complete proof of the one in \cite[page~22]{jones_arborealsurvey}:

\begin{lemma}[{see {\cite[page 22]{jones_arborealsurvey}}}]
\label{lemma_Jonesstrategy}
Let $G\le \Aut(T)$ be a closed subgroup. Assume that:
\begin{enumerate}[\normalfont(i)]
\item the fixed-point process $\{X_n\}_{n\ge 1}$ of $G$ is a martingale;
\item for any $r>0$ there exists $\epsilon:=\epsilon(r)$ and $m:=m(r)$ such that for infinitely many $n\ge 1$ we have $\mu(X_{n+m} = r \mid  X_n = r) \leq 1 - \epsilon$.
\end{enumerate}
Then, the fixed-point process of $G$ is eventually zero for almost all $g\in G$.
\end{lemma}
\begin{proof}
By assumption (i), the fixed-point process of $G$ is eventually constant for some constant $r:=r(g)$ for almost all $g\in G$. 

Let $\Lambda:=\{g\in G : X_n(g) \text{ becomes eventually constant}\}$ and for every $r\ge 0$ we further define $\Lambda_r:=\{g\in G: X_n(g) \text{ becomes eventually constant }r\}$. It is clear that $\Lambda:=\bigsqcup_{r\ge 0} \Lambda_r$. We also define for every $r\ge 0$ and any $n\ge 1$ the measurable subset $E_{n}^r:=\{g\in \Lambda: X_n(g)=r\}$. Then one has $\Lambda_r=\liminf E_n^r$.

We claim that for any $r> 0$ and any $m\ge 1$, we have 
$$\lim_{n\to{\infty}}\mu(E_{n}^r)=\lim_{n\to{\infty}}\mu(E_{n+m}^r\cap E_{n}^r)=\mu(\Lambda_r).$$

Therefore for $r>0$, if $\mu(\Lambda_r)>0$ we get
$$\lim_{n\to \infty} \mu(E_{n+m}^r\mid E_{n}^r)=\lim_{n\to \infty} \frac{\mu(E_{n+m}^r\cap E_{n}^r)}{\mu(E_n^r)}=\frac{\mu(\Lambda_r)}{\mu(\Lambda_r)}=1,$$
which contradicts assumption (ii). Thus $\mu(\Lambda_r)=0$ for every $r>0$ and the result follows as then $\mu(\Lambda_0)=\mu(\Lambda)=1$.

Thus, let us prove our claim. We prove it for $E_n^r$, as a similar argument also works for $E_{n+m}^r\cap E_n^r$. We only need to prove that $\lim_{n\to\infty} \mu(\Lambda_r\triangle E_n^r)=0$. Indeed,
\begin{align*}
\abs{\mu(\Lambda_r)-\mu(E_n^r)} &= \abs{\mu(\Lambda_r\setminus E_n^r) - \mu(E_n^r\setminus \Lambda_r)} \\
&\le \mu(\Lambda_r\setminus E_n^r)+\mu(E_n^r\setminus \Lambda_r)\\
&=\mu(\Lambda_r\triangle E_n^r).
\end{align*}

Now, observe that
$$\Lambda_r\triangle E_n^r = (\Lambda_r\setminus E_n^r)\sqcup (E_n^r\setminus \Lambda_r)\subseteq \bigsqcup_{r'\ge 0} \Lambda_{r'} \setminus E_n^{r'}.$$

Since the measure $\mu$ is finite and $\Lambda=\bigsqcup_{r'\ge 0} \Lambda_{r'}$, by dominated convergence theorem, we may interchange the summation and the limit to obtain
$$\lim_{n\to\infty} \mu(\Lambda_r\triangle E_n^r) \leq \lim_{n\to\infty} \sum_{r'\ge 0} \mu(\Lambda_{r'}\setminus E_n^{r'})=\sum_{r'\ge 0} \lim_{n\to\infty} \mu(\Lambda_{r'}\setminus E_n^{r'}).$$
Hence, it is enough to show that $\lim_{n\to\infty} \mu(\Lambda_{r'}\setminus E_n^{r'})=0$ for every $r'\ge 0$. Let us write $F_k^{r'}:=\bigcap_{\ell\ge k} E_\ell^{r'}$ for $k\ge 1$ and $F_{0}^{r'}=\emptyset$. Then 
$$\Lambda_{r'}=\liminf_{k\to \infty} E_k^{r'}=\bigcup_{k\ge 1}F_k^{r'}=\bigsqcup_{k\ge 1} F_{k}^{r'}\setminus F_{k-1}^{r'}=F_{n}^{r'}\sqcup \left(\bigsqcup_{k\ge n+1} F_{k}^{r'}\setminus F_{k-1}^{r'}\right).$$
Since $\mu(\Lambda_{r'})$ is finite we get that $\lim_{n\to \infty}\mu(\bigsqcup_{k\ge n+1} F_{k}^{r'}\setminus F_{k-1}^{r'})=0$ and therefore
\begin{align*}
0\le \lim_{n\to\infty} \mu(\Lambda_{r'}\setminus E_n^{r'})\le \lim_{n\to\infty} \mu(\Lambda_{r'}\setminus F_n^{r'})=\lim_{n\to \infty}\mu(\bigsqcup_{k\ge n+1} F_{k}^{r'}\setminus F_{k-1}^{r'})=0,
\end{align*}
as $F_{n}^{r'}\subseteq E_n^{r'}$.
\end{proof}

\subsection{Proof of \cref{theorem: FPP SSF}}

We shall see that super strongly fractal groups satisfy both assumptions in \cref{lemma_Jonesstrategy}. To check assumption (i) in \cref{lemma_Jonesstrategy}, we shall use the following:

\begin{theorem}[{see {\cite[Theorem 5.7]{BridyJones2022} and \cite{Jones_density_prime_divisors}}}]
\label{teo_condition_martingale}
Let $G \leq \Aut(T)$ be a closed subgroup. Then, $\set{X_n}_{n \geq 0}$ is a martingale if and only if for all $n \geq 1$, the action of $\St_G(n-1)$ on $ T_v^1$ is transitive for any $v \in \mathcal{L}_{n-1}$.
\end{theorem}

\begin{lemma}
If $G$ is super strongly fractal, then the fixed point process $\set{X_n}_{n \ge 1}$ is a martingale.
\label{lemma_ssf_is_martingale}
\end{lemma}

\begin{proof}
For every $n\ge 1$ and every vertex $v\in \mathcal{L}_{n-1}$ we have $\varphi_v(\mathrm{St}_G(n-1))=G$ by the super strong fractality of $G$. Thus, the result follows by \cref{teo_condition_martingale} as $G$ acts transitively on the first level of the tree.
\end{proof}

For any $n\ge 1$ and $A\subseteq \pi_n(G)$ we define the \textit{cone} $C_{A}\subseteq G$ via $C_A:=\pi_n^{-1}(A)$. The proof of \cref{theorem: FPP SSF} relies on the following key lemma proved by the first author in \cite{JorgeCyclicity}:

\begin{lemma}[{see {\cite[Lemma 4.2]{JorgeCyclicity}}}]
\label{lemma_measure_of_intersection}
    Let $G\le \Aut(T)$ be super strongly fractal and let us fix $m,n\ge 1$. Then, for any pair $A\subseteq \pi_n(G)$ and $B\subseteq \pi_m(G)$ and any $v\in \mathcal{L}_n$ we have
    $$\mu(C_A\cap \varphi_v^{-1}(C_B))= \mu(C_A)\cdot\mu(C_B).$$
\end{lemma}

\cref{lemma_measure_of_intersection} implies that the measure-preserving dynamical system associated to a super strongly fractal group is mixing. This emphasizes the key role that the ergodic theory developed in \cite{JorgeCyclicity} plays in the fixed-point proportion of self-similar groups.

\begin{proof}[Proof of \cref{theorem: FPP SSF}]
To prove \cref{theorem: FPP SSF} it is enough to show that the fixed-point process of $G$ is eventually zero for almost all $g\in G$. Thus, we just need to prove that~$G$ satisfies condition (ii) in \cref{lemma_Jonesstrategy}, as condition (i) in \cref{lemma_Jonesstrategy} was already proved in \cref{lemma_ssf_is_martingale}. Fix $r>0$ and $m:=\lceil \log_d(r)\rceil+1$, which only depends on $r$. Now, for any $n\ge 1$ consider the subset of $n$-patterns
$$A_n:=\set{g \in \pi_n(G): g \text{ has exactly $r$ fixed points in } \mathcal{L}_n }.$$
Note that for any $a\in A_n$ there exists a vertex $v_a\in \mathcal{L}_n$ fixed by $a$. By \cref{lemma_measure_of_intersection} we obtain that
$$\mu(C_{\{a\}}\cap \varphi_{v_a}^{-1}(C_{\{\mathrm{id}\}}))=\mu(C_{\{a\}})\cdot \mu(C_{\{\mathrm{id}\}}).$$
Recall also that
$$\mu(C_{A_n})=\sum_{a\in A_n} \mu(C_{\{a\}})$$
as $C_{A_n}=\bigsqcup_{a\in A_n}C_{\{a\}}$ is a disjoint union of measurable subsets.

Now, for every $a\in A_n$, an element $g\in C_{\{a\}}\cap \varphi_{v_a}^{-1}(C_{\{\mathrm{id}\}})$ satisfies both that $X_n(g)=r$ and $X_{n+m}(g)>r$ by our choice of $m$, as $g$ fixes at least $d^m>r$ vertices at level $n+m$. Therefore
\begin{align*}
\mu(X_{n+m} = r \mid  X_n = r)&\le 1-\mu(X_{n+m}>r\mid X_n=r)\\
&\le 1- \frac{\mu\big(\bigsqcup_{a\in A_n}(C_{\{a\}}\cap \varphi_{v_a}^{-1}(C_{\{\mathrm{id}\}}))\big)}{\mu(C_{A_n})}\\
&=1- \frac{\sum_{a\in A_n}\mu(C_{\{a\}}\cap \varphi_{v_a}^{-1}(C_{\{\mathrm{id}\}}))}{\mu(C_{A_n})}\\
&= 1- \frac{\sum_{a\in A_n}\mu(C_{\{a\}})\cdot \mu(C_{\{\mathrm{id}\}})}{\mu(C_{A_n})}\\
&= 1- \frac{\mu(C_{A_n})\cdot \mu(C_{\{\mathrm{id}\}})}{\mu(C_{A_n})}\\
&= 1-|\pi_m(G)|^{-1},
\end{align*}
and condition (ii) in \cref{lemma_Jonesstrategy} is satisfied by $G$ with $\epsilon:=|\pi_m(G)|^{-1}$.
\end{proof}

Many well-known groups studied in geometric group theory are super strongly fractal. For example, the first Grigorchuk group \cite[Proposition 5.3]{UriaAlbizuri2016}, the Basilica group \cite[Theorem 1]{Basilica1} and its generalization to the $p$-adic tree \cite[Theorem 4.5]{pBasilica2021}, and the GGS-groups with non-constant defining vector \cite[Proposition 5.1]{UriaAlbizuri2016} among others. Therefore \cref{theorem: FPP SSF} implies directly that all these well-known groups have zero fixed-point proportion.

\section{The fixed-point proportion of iterated monodromy groups of dynamically exceptional polynomials}
\label{sec_FPP_dyn_exc_poly}

In this section, we give a brief introduction to iterated monodromy groups and prove \cref{theorem: dynamically exceptional FPP=0}, namely, that there exist dynamically exceptional polynomials of each degree $d\ge 3$ such that their iterated monodromy groups have null fixed-point proportion. For further details on iterated monodromy groups the reader is encouraged to read \cite[Chapters 5 \& 6]{Self_similar_groups}.

\subsection{A primer on iterated monodromy groups}
\label{sec_primer_IMG}

Let $f: \PP^1(\CC) \rightarrow \PP^1(\CC)$ be a polynomial of degree $d\geq 2$. The \textit{local degree} of $f$ at $z$, denoted $e_f(z)$, is the smallest $k\ge 0$ such that the $k$-th derivative of $f$ is nonzero at $z$. Let $C_f := \set{z \in \PP^1(\CC): e_f(z) > 1}$ be the set of critical points. The set $C_f$ is finite as $f$ is a polynomial. Then, let us define $P_f := \set{f^n(c): c \in C_f \text{ and } n \geq 1}$ the \textit{post-critical set} of~$f$. We say that $f$ is \textit{post-critically finite} if the set $P_f$ is finite.

Note that $f^n: \PP^1(\CC) \setminus f^{-n}(P_f) \rightarrow \PP^1(\CC) \setminus P_f$ is an unramified covering for all $n \geq 1$. Indeed $f^n$ is ramified if there exists $z \in \PP^1(\CC) \setminus f^{-n}(P_f)$ such that $(f^n)'(z) = 0$, or, by the chain rule, if there exists $0\le i \le n-1$ such that $f^i(z) \in C_f$. However, this implies that $f^n(z) \in f^{n-i}(C_f) \subseteq P_f$, so $z \in f^{-n}(P_f)$.

Let us consider $t \in \PP^1(\CC) \setminus P_f$, a loop $\gamma$ starting and ending at $t$ and $\widetilde{t} \in f^{-n}(t)$. By \cite[Proposition 1.34]{Hatcher_Algebraic_topology}, there exists a unique lift $\widetilde{\gamma}$ of $\gamma$ starting at $\widetilde{t}$. Since $\gamma$ is a loop, the ending point of $\widetilde{\gamma}$ is necessarily another point in $f^{-n}(t)$. Moreover, if $\gamma'$ is another loop starting and ending at $t$ that is homotopically equivalent to $\gamma$, the ending point of $\widetilde{\gamma'}$, the unique lift of $\gamma'$, has the same endpoint as $\widetilde{\gamma}$. Thus, we have well-defined actions of $\pi_1(\PP^1(\CC) \setminus P_f, t)$ on $f^{-n}(t)$ for each natural number $n \geq 1$. 

These actions are coherent in the sense that if $\widetilde{\gamma}_n$ is a lift of $\gamma$ starting at $\widetilde{t}_n \in f^{-n}(t)$ and $\widetilde{\gamma}_{n+1}$ is a lift of $\gamma$ starting at $\widetilde{t}_{n+1} \in f^{-1}(\widetilde{t}_n)$ then by uniqueness of the lifts, we have $f(\widetilde{\gamma}_{n+1}) = \widetilde{\gamma}_n$ and therefore the ending point of $\widetilde{\gamma}_{n+1}$ is mapped via~$f$ to the ending point of $\widetilde{\gamma}_n$. 

The action of $\pi_1(\PP^1(\CC) \setminus P_f, t)$ on the disjoint union $\bigsqcup_{n \geq 0} f^{-n}(t)$ is called the \textit{monodromy action}.

As $t \notin P_f$, the set $f^{-n}(t)$ has cardinality $d^n$ for each $n\ge 1$ and we may identify the set of all preimages $\bigsqcup_{n\ge 0} f^{-n}(t)$ with the $d$-regular tree $T$. Consequently, we obtain an action of $\pi_1(\PP^1(\CC) \setminus P_f, t)$ on $T$. However, this action is not necessarily faithful. Thus, we may quotient $\pi_1(\PP^1(\CC) \setminus P_f, t)$ with the kernel of this action to obtain a subgroup of the automorphism group of the $d$-regular tree $T$. This group is denoted $\IMG(f)$ and called the \textit{iterated monodromy group} of $f$.

Note that if we restrict our attention to post-critically finite polynomials, then $\PP^1(\CC) \setminus P_f$ is path-connected and the induced action of $\pi_1(\PP^1(\CC) \setminus P_f, t)$ on $T$ is unique up to conjugation in $\Aut(T)$ for any choice of the basepoint $t$. 

\subsection{Hyperbolic polynomials}

Let $f$ be a post-critically finite complex polynomial. Following \cite{DoaudyHubbard1993}, we define its \textit{Thurston orbifold} as $(\mathbb{P}^1(\mathbb{C}),\nu_f)$ for the map $\nu_f:\mathbb{P}^1(\mathbb{C})\to \mathbb{N}\cup \{\infty\}$, where
$$\nu_f(z):=\mathrm{lcm}\{\nu_f(w)e_f(w) : w\in f^{-n}(z) \text{ for some }n\ge 1\},$$
unless $z$ belongs to a super-attracting periodic orbit (i.e. a periodic orbit containing a critical point of $f$), in which case $\nu_f(z)=\infty$. Note that $\nu_f(z)=1$ for all $z\notin P_f$. We say that $f$ is \textit{hyperbolic} if 
$$\chi_f := 2- \sum_{z \in P_f} \left(1 - \frac{1}{\nu_f(z)} \right) < 0.$$



\subsection{Construction of the iterated monodromy group of a polynomial using spiders}
\label{section: iterated monodromy group}

We shall follow the approach in \cite[Section 6.8]{Self_similar_groups} to compute the iterated monodromy group of a given post-critically finite polynomial $f: \PP^1(\CC) \rightarrow \PP^1(\CC)$ using spiders. Let $d:=\mathrm{deg}(f)$. A \textit{spider} is a collection of simple curves $\set{\gamma_p: p \in P_f}$, where for each $p\in P_f$ the curve $\gamma_p: [0,1] \rightarrow \big(\PP^1(\CC) \setminus (C_f \cup P_f) \big)\bigcup  \set{p,\infty}$ starts at $p$ and ends at $\infty$. If $c \in C_f$ is a critical point with local degree $e_f(c)$, then $f^{-1}(\gamma_{f(c)})$ consists of $e_f(c)$ distinct curves starting at $c$ and ending at $\infty$. The curves in $\set{f^{-1}(\gamma_{f(c)}): c \in C_f}$ cut the plane $\CC$ in $d$ \textit{sectors} $S_1,\dotsc, S_d$.

We fix the alphabet $X:=\set{1,\dotsc, d}$ and define a group $G_f$. The generators of~$G_f$ will be $\set{g_p}_{p \in P_f}$. To define their actions, take $p \in P_f$ and let $\alpha_p$ be a small loop around $p$ in the positive direction. Let $y \in f^{-1}(p)$ and let also $\alpha_{p,y}$ be one of the $f$-preimages of $\alpha_p$ corresponding to a small loop around $y$.

\begin{enumerate}[\normalfont(i)]
    \item \textit{The point $y$ is not a critical point}: Then its local degree is $1$ and there exists $x \in \set{1,\dotsc,d}$ such that $y$ belongs to the sector $S_x$. In this case $g_p(x) = x$ and 
\begin{align*}
g_p|_x = \left \{ \begin{matrix} 
g_y & \mbox{if $y \in P_f$,}  \\
1 & \mbox{if $y \notin P_f$.}
\end{matrix}\right.
\end{align*}

\item \textit{The point $y$ is a critical point but not post-critical}: Then $y$ is in the boundary of $d':=e_f(y)$ sectors $S_{x_1},\dotsc,S_{x_{d'}}$. Numbering $x_1,\dotsc,x_{d'}$ according to the order in which $\alpha_{p,y}$ meets them, let us define $g_p(x_i) = x_{i+1}$ for $1 \leq i \leq d'-1$ and $g_p(x_{d'}) = x_1$, and $g_p|_{x_i} = 1$ for $1\le i\le d'$.

\item \textit{The point $y$ is both a critical and post-critical point}: Then $y$ is in the boundary of $d':=e_f(y)$ sectors and there is also a path $\gamma_y$ in the spider being adjacent to two of those sectors. Let us number again $x_1,\dotsc,x_{d'}$ according to the order in which $\alpha_{p,y}$ meets them. If $\gamma_p$ is adjacent to $S_{x_1'}$ and $S_{x_d'}$ then $g_p(x_i) = x_{i+1}$ for $1 \leq i \leq d'-1$ and $g_p(x_{d'}) = x_1$. Furthermore
\begin{align*}
g_p|_{x_i} = \left \{ \begin{matrix} 
g_y & \mbox{if $i = d'$,}  \\
1 & \mbox{if $i \neq d'$.}
\end{matrix}\right.
\end{align*}
\end{enumerate}

Note that for hyperbolic post-critically finite polynomials the group~$G_f$ coincides with the corresponding iterated monodromy group:

\begin{theorem}[{see {\cite[Theorems 6.8.3 and 6.9.1]{Self_similar_groups}}}]
Let $f: \PP^1(\CC) \rightarrow \PP^1(\CC)$ be a hyperbolic post-critically finite polynomial. Then $\IMG(f) = G_f$.
\label{theorem: IMG = Af hyperbolic map}
\end{theorem}

\subsection{Dynamically exceptional polynomials}
\label{section: dynamically exceptional polynomials}

Recall that a post-critically finite complex polynomial $f$ of degree $d \geq 2$ is \textit{dynamically exceptional} if there exists a non-empty set $E \subseteq \PP^1(\CC)$ such that $f^{-1}(E) \setminus C_f = E$. Notice that this condition implies that $E$ does not contain critical points. As it was mentioned in \cref{section: introduction}, the set $E$ can have one or two elements. The cases of dynamically exceptional polynomials with $\# E = 2$ were completely solved by Jones in \cite{jones2012fixedpointfree}. If $E = \{z_0\}$, then $f^{-1}(z_0)$ consists of $z_0$ and a subset of $C_f$, but $z_0 \notin C_f$. Conjugating $f$ by a translation, we may assume that $z_0 = 0$ and $f$ has the form 
$$f(z) = z(z-a_1)^{k_1}\dotsb(z-a_m)^{k_m}$$ 
with $a_1,\dotsc,a_m \in \CC \setminus \set{0}$ and $k_i \geq 2$ for all $i = 1,\dotsc,m$.

In the above setting, i.e. when $f$ is dynamically exceptional with $\# E = 1$, no previous result is known to the best of our knowledge.

\subsection{Proof of \cref{theorem: dynamically exceptional FPP=0}}

We shall use the following lemma:

\begin{lemma}[{see {\cite[Corollary 6.7.7 and Proposition 6.8.2]{Self_similar_groups}}}]
\label{lemma: cycles}
    Let $f:\mathbb{P}_1(\mathbb{C})\to \mathbb{P}_1(\mathbb{C})$ be a post-critically finite complex polynomial of degree $d$ and let $g_1,\dotsc,g_\ell$ be the generators of the group $G_f$. Then, a product of all the generators $g_1,\dotsc,g_\ell$ in any order is a level-transitive element in the automorphism group of the $d$-regular tree. In particular $G_f$ is level-transitive.
\end{lemma}

We are now ready to prove \cref{theorem: dynamically exceptional FPP=0}, namely that for each $d \geq 3$, there exists a dynamically exceptional polynomial $f$ whose degree is $d$ and $\FPP(\IMG(f)) = 0$.

\begin{proof}[Proof of \cref{theorem: dynamically exceptional FPP=0}]
Let us fix $d\ge 3$ and consider $f(z) = z(z-a)^{d-1}$ where $a \neq 0$. The critical points of $f$, i.e. when $f'(z) = 0$, are $a$ and $a/d$. The local degree of~$a$ is clearly $d-1$. In order to calculate the local degree of $a/d$ we may use the Riemann-Hurwitz formula \cite[Theorem 1.1]{Silverman2007}. Note that $\infty$ is a totally ramified point as $f$ is a polynomial, so $e_f(\infty)=d$. Then we get the equality
$$2d-2 = (e_f(a)-1) + (e_f(\infty) -1) + (e_f(a/d)-1),$$ 
which yields that $e_f(a/d) = 2$. 

The point $0$ is a fixed point for $f$. Thus $a$ is strictly preperiodic with period $1$, as $a \mapsto 0$. Now, let us assume that $a$ is chosen such that $a/d$ is a fixed point , i.e. $f(a/d) = a/d$. Solving for $a$, we get 
$$a := \zeta \frac{d}{(1-d)},$$
where $\zeta$ is any $(d-1)$st root of unity. In particular $a \neq a/d$. We now use the algorithm presented in \cref{section: iterated monodromy group} to calculate $G_f$.

As $P_f = \set{0,a/d}$, the group $G_f$ is generated by two elements $g:=g_0$ and $h:=g_{a/d}$, corresponding to the post-critical points 0 and $a/d$ respectively. In order to calculate the sections of $h$, we observe that the preimages of $a/d$ are: $a/d$ with multiplicity~2, as $a/d$ is a critical fixed-point with local degree $2$; and $d-2$ different complex points, since the only critical points of $f$ are $a$ and $a/d$ and $a$ is not mapped to $a/d$. Hence~$h$ has $d-1$ sections equal to $1$ and one section equal to $h$. Moreover, the action of $h$ on the first level is a transposition $(i_1\,i_2)$, where exactly one section among $h|_{i_1}$ and $h|_{i_2}$ is equal to $h$. By relabelling vertices if needed, we may assume that
$$h=(h,1,\dotsc,1)(1\, 2).$$

To calculate the sections of $g$, we observe that the preimages of $0$ via $f$ are $0$ with multiplicity 1 and $a$ with multiplicity $d-1$. Hence, all the sections of $g$ are trivial except for one which is $g$. Moreover, since $a \notin P_f$, the automorphism $g$ acts on the first level as a $d-1$ cycle permuting the vertices at which it has a trivial section and fixing the vertex $v$ where $g|_v=g$. By \cref{lemma: cycles}, we must have $v\in \{1,2\}$, since~$h$ acts on the first level as $(1\, 2)$. In fact, we may further assume that $v=1$ by considering the relabelling $1\leftrightarrow 2$ and considering $h^{-1}$ instead if $v=2$. In other words 
$$g=(g,1,\dotsc,1)(2\,\dotsb \, d).$$

We claim that $G_f = \langle g,h\rangle$ is super strongly fractal. As $G_f$ acts transitively on the first level of $T$ by \cref{lemma: cycles}, we only need to prove that $\varphi_v(\mathrm{St}_{G_f}(n))=G_f$ holds for every vertex $v$ at the $n$th level for each $n\ge 1$. Now observe that $h^2\in \mathrm{St}_{G_f}(1)$ and it has precisely two sections equal to $h$ and the rest of the sections are trivial. Hence by induction, $h^{2^n}\in \mathrm{St}_{G_f}(n)$ and it has precisely $2^n$ sections at the $n$th level equal to $h$ and the rest are trivial. Therefore $h \in \varphi_v(\mathrm{St}_{G_f}(n))$ for any $v$ among these $2^n$ vertices $\{\widetilde{v}_1,\dotsc,\widetilde{v}_{2^n}\}$ at the $n$th level.

Note now that $(gh)^d,(hg)^d\in \mathrm{St}_{G_f}(1)$ and every section of these elements at the first level is either $gh$ or $hg$. Hence, arguing by induction again, we obtain that $(gh)^{d^n},(hg)^{d^n}\in \mathrm{St}_{G_f}(n)$ and their sections at the $n$th level are either $gh$ or $hg$. Thus, either $gh\in \varphi_w(\mathrm{St}_{G_f}(n))$ or $hg\in \varphi_w(\mathrm{St}_{G_f}(n))$ for any~$w$ at the $n$th level and so in particular we obtain that $\varphi_v(\mathrm{St}_{G_f}(n))=G_f$ for $v\in \{\widetilde{v}_1,\dotsc,\widetilde{v}_{2^n}\}$. Then, the condition $\varphi_w(\mathrm{St}_{G_f}(n))=G_f$ holds for every vertex $w$ at the $n$th level for each $n\ge 1$ and the group $G_f$ is super strongly fractal.

We finally need to prove that $G_f = \IMG(f)$. By \cref{theorem: IMG = Af hyperbolic map}, we only need to prove that $f$ is hyperbolic. We calculate $\nu_f$ first. The function $\nu_f$ may only be different to $1$ at the points $0$, $\infty$ and $a/d$. The preimages of $0$ are $\set{0,a}$ \and $e_f(0)=1$ and $e_f(a)=d-1$. Thus $\nu_f(0)=d-1$. As both $a/d$ and $\infty$ are critical points in a periodic orbit, we get $\nu_f(a/d)=\nu_f(\infty)=\infty$. Therefore 
\begin{align*}
\chi_f &=2-\sum_{z\in P_f}\left(1-\frac{1}{\nu_f(z)}\right)\\
&=2-\left(1-\frac{1}{\nu_f(0)}\right)-\left(1-\frac{1}{\nu_f(\infty)}\right)-\left(1-\frac{1}{\nu_f(a/d)}\right)\\
&= 2 - \left( 1 - \frac{1}{d-1} \right) -1 - 1 \\
&=  - \left( 1 - \frac{1}{d-1} \right)<0,
\end{align*}
and consequently $f$ is hyperbolic. 
\end{proof}


\section{Examples}
\label{section: Scope of theorem}

First, we show that \cref{theorem: FPP SSF} cannot be extended neither to fractal groups nor to groups of finite type. In fact, we show that branching properties of self-similar groups are not related to the fixed-point proportion of the group in general. 

\subsection{Fractal groups}

In this example we show an example of a fractal group whose fixed-point proportion is positive. Indeed, in the binary tree, consider $G := \<a,b>$ where $a = (1,1)\sigma$ and $b = (b,a)$ being $\sigma$ the non-trivial permutation in $\Sym(2)$. 

Then $b, aba \in \St_G(1)$ and if $v$ is the leftmost vertex at the first level then $b|_v = b$ and $(aba)|_v = a$. Thus $G$ is strongly fractal as it is also level-transitive. However, the group $G$ corresponds to the iterated monodromy group of the Chebyshev polynomial $T_2 = 2x^2-1$ by \cite[Proposition 6.12.6]{Self_similar_groups}, whose fixed-point proportion is $1/4$ by \cite[Proposition 1.2]{jones2012fixedpointfree}. 

\subsection{Groups of finite type}
 
A group $G \leq \Aut(T)$ is said to be of \textit{finite type} if there exists $D \ge 1$, called the \textit{depth}, and $\mathcal{P} \leq \Aut(T^D)$ such that $G = G_\mathcal{P}$ where $$G_\mathcal{P} := \set{g \in \Aut(T): g|_v^D \in \mathcal{P}\text{ for all }v\in T}.$$ 
Groups of finite type are self-similar and closed; see for example \cite[page 3]{Onfinitegeneration}. In particular, if $D = 1$, then $g|_v^1 \in \mathcal{P}$ at any vertex $v \in T$ and the resulting group is the \textit{iterated wreath product} $W_\mathcal{P}$. By \cite[Corollary 2.7]{AbertVirag2004}, level-transitive iterated wreath products have null fixed-point proportion. Note that this also follows from \cref{theorem: FPP SSF} as level-transitive iterated wreath products are super strongly fractal. Groups of finite type represent the most immediate generalization of iterated wreath products.

However, there are examples of groups of finite type whose fixed-point proportion is not zero, as proved by the second author in \cite{Radi2024}. For instance, in the $3$-adic tree, consider $$G_{\mathrm{Alt}(3)}^{\Sym(3)} := \set{g \in \Aut(T): (g|_v^1)(g|_w^1)^{-1} \in \mathrm{Alt}(3) \text{ for all } v,w \in T},$$ where $\mathrm{Alt}(3)\le \mathrm{Sym}(3)$ is the alternating group. The elements in the group $G_{\mathrm{Alt}(3)}^{\Sym(3)}$ are those whose labels all lie on the same coset of $\mathrm{Alt}(3)$ in $\mathrm{Sym}(3)$. Since the alternating group $\mathrm{Alt}(3)$ is transitive, the group $G_{\mathrm{Alt}(3)}^{\Sym(3)}$ is level-transitive and its fixed-point process is a martingale. By \cite[Theorem 2]{Radi2024} the group $G$ has fixed-point proportion $1/2$. Note that the group $G_{\mathrm{Alt}(3)}^{\Sym(3)}$ is fractal; see \cite[Proposition 2.3]{JorgeFPP} for a proof.

\bibliographystyle{unsrt}

\begin{thebibliography}{1}

\bibitem{AbertVirag2004}
M. Abért, and B. Virág, Dimension and randomness in groups acting on rooted trees, \textit{J. Amer. Math. Soc.} \textbf{18(1)} (2005), 157--192.

\bibitem{Onfinitegeneration}
I. Bondarenko, and I. Samoilovych, On finite generation of self-similar groups of finite type, \textit{Internat. J. Algebra Comput.} \textbf{23(1)} (2013), 69--79.

\bibitem{BridyJones2022}
A. Bridy, R. Jones, G. Kelsey, and R. Lodge, Iterated monodromy groups of rational functions and periodic points over finite fields, \textit{Math. Ann.} \textbf{390(1)} (2024), 439--475.

\bibitem{pBasilica2021}
E. Di Domenico and G.\,A. Fern\'andez-Alcober and M.\,L. Noce and A. Thillaisundaram, $p$-Basilica groups, \textit{Mediterr. J. Math.} \textbf{19(6)} (2022), 275.

\bibitem{DoaudyHubbard1993}
A. Douady and J.\,H. Hubbard, A proof of {T}hurston's topological characterization of rational functions, \textit{Acta Math.} \textbf{171(2)} (1993), 263--297.

\bibitem{JorgeFPP}
J. Fariña-Asategui, Arboreal Galois representations of rational functions: fixed-point proportion and the extension problem, arXiv preprint: 2601.19414.

\bibitem{JorgeCyclicity}
J. Fariña-Asategui, Cyclicity, hypercyclicity and randomness in self-similar groups, \textit{Monatsh. Math.} \textbf{207} (2025), 275--292.

\bibitem{Farina2024Restricted}
J. Fariña-Asategui, Restricted Hausdorff spectra of $q$-adic automorphisms, \textit{Adv. Math.} \textbf{472} (2025), 110294.

\bibitem{Corrigendum}
J. Fariña-Asategui, R. Jones and S. Radi, Corrigendum and addendum to ``Fixed-point-free elements of iterated monodromy groups" [Trans. of the American Mathematical Society Vol. 367, No. 3 (2015), 2023-2049], arXiv preprint: 2602.18654.

\bibitem{Fried1986}
M.\,D. Fried, Field Arithmetic, \textit{Springer Berlin Heidelberg}, 11 (2023), xxxi--827.

\bibitem{Grig_1980}
R.\,I. Grigorchuk, On Burnside’s problem for periodic groups, \textit{Funktsional. Anal. i Prilozhen} \textbf{14(1)} (1980), 53--54.

\bibitem{Liftability}
R.\,I. Grigorchuk and D. Savchuk, Liftable self-similar groups and scale groups, \textit{Trans. Amer. Math. Soc.}, to appear.

\bibitem{Basilica1}
R.\,I. Grigorchuk and A. \.{Z}uk, Spectral properties of a torsion-free weakly branch group defined by a three state automaton, \textit{Amer. Math. Soc. Providence, RI}, \textbf{298} (2002), 57--82.

\bibitem{Grimmett_Random_processes}
G.\,R. Grimmett and D.\,R. Stirzaker, Probability and random processes, \textit{Oxford University Press, Oxford}, 60--01 (2020), xii--669.

\bibitem{GGS1983}
N. Gupta and S. Sidki, On the Burnside Problem for Periodic Groups, \textit{Math. Z.} \textbf{182} (1983), 385--388.

\bibitem{jones2012fixedpointfree}
R. Jones, Fixed-point-free elements of iterated monodromy groups, \textit{Trans. Amer. Math. Soc.} \textbf{367(3)} (2015), 2023--2049.

\bibitem{jones_arborealsurvey}
R. Jones, Galois representations from pre-image trees: an arboreal survey, \textit{Pub. Math. Besancon} (2013), 107--136.

\bibitem{JonesManes2012}
R. Jones and M. Manes, Galois theory of quadratic rational functions, \textit{Comment. Math. Helv.} \textbf{89(1)} (2014), 173--213.

\bibitem{Jones2007}
R. Jones, Iterated Galois towers, their associated martingales, and the $p$-adic Mandelbrot set, \textit{Compos. Math.} \textbf{143(5)} (2007), 1108--1126. 

\bibitem{Jones_density_prime_divisors}
R. Jones, The density of prime divisors in the arithmetic dynamics of quadratic polynomials, \textit{J. London Math. Soc.} \textbf{78(2)} (2008), 523--544.

\bibitem{juul2014wreath}
J. Juul, P. Kurlberg, K. Madhu, and T.\,J. Tucker, Wreath products and proportions of periodic points, \textit{Int. Math. Res. Not. IMRN} \textbf{13} (2016), 3944--3969.

\bibitem{Hatcher_Algebraic_topology}
A. Hatcher, Algebraic topology, \textit{Cambridge University Press}, (2002), xii--544.

\bibitem{MakarovSmirnov1996}
N. Makarov and S. Smirnov, On “thermodynamics” of rational maps I. Negative spectrum, \textit{Commun. Math. Phys.} \textbf{211(3)} (2000), 705--743.

\bibitem{Self_similar_groups}
V. Nekrashevych, Self-similar groups, \textit{American Mathematical Society, Providence, RI} 117 (2005), xii--231.

\bibitem{Odoni_1985}
R.\,W.\,K. Odoni, On the Prime Divisors of the Sequence $w_{n+1} = 1 + w_1 \dots w_n$, \textit{J. London Math. Soc.} \textbf{32(2)} (1985), 1--11.

\bibitem{GGS_contracting}
J.\,M. Petschick, On conjugacy of GGS-groups, \textit{J. Group Theory} \textbf{22(3)} (2019), 347--358.

\bibitem{Radi2024}
S. Radi, Branch iterated Galois groups with positive fixed-point proportion and positive Hausdorff dimension, \textit{} arXiv preprint: 2602.13428.

\bibitem{Ribes2000}
L. Ribes, P. Zalesskii, Profinite Groups, \textit{Springer-Verlag, Berlin}, 40 (2010), xvi--464.

\bibitem{Silverman2007}
J.\,H. Silverman, The arithmetic of dynamical systems, \textit{Springer, New York}, 241 (2007), x--511.

\bibitem{UriaAlbizuri2016}
J. Uria-Albizuri, On the concept of fractality for groups of automorphisms of a regular rooted tree, \textit{Reports@SCM} \textbf{2} (2016), 33--44.

\end{thebibliography}

\end{document}